\documentclass{amsart}
\usepackage{graphics}
\newtheorem{theorem}{Theorem}[section]
\newtheorem{lem}[theorem]{Lemma}
\newtheorem{prop}[theorem]{Proposition}
\newtheorem{cor}[theorem]{Corollary}
\theoremstyle{definition}

\newtheorem{example}[theorem]{Example}

\newtheorem{conj}[theorem]{Conjecture}

\theoremstyle{remark}

\numberwithin{equation}{section}

\newcommand{\half}{\frac{1}{2}}
\newcommand{\Hor}{{\mathcal{H}}}
\newcommand{\V}{{\mathcal{V}}}
\newcommand{\ra}{\rightarrow}
\newcommand{\Soul}{\Sigma}

\newcommand{\lb}{\langle}
\newcommand{\rb}{\rangle}

\newcommand{\SO}{\text{SO}}
\newcommand{\so}{\text{so}}

\newcommand{\R}{\mathbb{R}}
\newcommand{\D}{\mathcal{D}}

\newcommand{\RN}{R^{\nabla}}

\newcommand{\ON}{\overline{\nabla}}
\newcommand{\BX}{\overline{X}}
\newcommand{\BY}{\overline{Y}}
\newcommand{\bp}{\overline{p}}

\newcommand{\nbp}{\nu_p(\Sigma)}
\newcommand{\nb}{\nu(\Sigma)}
\newcommand{\U}{\mathcal{U}}

\begin{document}

\title{Rigidity for Nonnegatively curved metrics on $S^2\times\R^3$}
\author{Kristopher Tapp}
\address{Department of Mathematics\\ Bryn Mawr College\\
         Bryn Mawr, PA 19010}
\email{ktapp@brynmawr.edu}

\subjclass{Primary 53C20}
\date{\today}
\keywords{nonnegative curvature, soul}

\begin{abstract}
We address the question: how large is the family of complete metrics with nonnegative sectional curvature on
$S^2\times\R^3$?  We classify the connection metrics, and give several examples of non-connection metrics.
We provide evidence that the family is small by proving some rigidity results for metrics more general than
connection metrics.
\end{abstract} \date{\today} \maketitle


\section{Introduction}

According to~\cite{GT},
the space $S^2\times\R^2$ with an arbitrary complete metric of nonnegative sectional curvature
must be a product metric or be isometric to a Riemannian quotient of the form $$((S^2,g_0)\times(\R^2,g_f)\times\R)/\R,$$
were $g_0$ and $g_f$ are $\R$-invariant metrics.  These metrics are
``rigid at the soul'', meaning the following inequality from~\cite{T} becomes an equality in
the case when $M$ is diffeomorphic to $S^2\times\R^2$:
\begin{prop}\label{P:condition}\cite{T}
If $M$ is an open manifold of nonnegative curvature with soul $\Sigma$, then for
any $p\in\Sigma$, orthonormal vectors $X,Y\in T_p\Sigma$, and orthonormal $W,V\in\nbp$ = the normal space to $T_p\Sigma$ in $T_pM$,
$$
\lb (D_XR)(X,Y)W,V\rb^2\leq(|R(W,V)X|^2+(2/3)(D_{X,X}R)(W,V,V,W))\cdot R(X,Y,Y,X).
$$
\end{prop}

In general, the metric sphere, $S_{\epsilon}$, of small radius $\epsilon$ about a soul inherits a metric of
nonnegative curvature~\cite{GW}.  If the above inequality is strict, then $S_{\epsilon}$ is positively curved~\cite{T}.
If there is a single $(p,W)\in\nb$ such that the inequality is strict for all
$X,Y,V$ at $p$, then $\epsilon W$ is a point of $S_{\epsilon}$ with positive curvature, so $S_{\epsilon}$ is quasi-positively
curved, which means that it has nonnegative curvature, and positive curvature at a point
(in this case, we will say that the inequality is ``quasi-strict'').

The question whether $S^2\times S^2$ admits
positive or quasi-positive curvature is a long-standing unsolved
problem in Riemannian geometry, which motivates our study of metrics of
nonnegative curvature on $S^2\times\R^3$.  Although
we do not completely classify such metrics, we demonstrate special cases under which the inequality rigidly
determines what the metrics may look like at the soul, which provides evidence that the family of metrics is small.

We first classify the connection metrics, i.e., metrics with totally geodesic Sharafutdinov fibers.
It is convenient to simultaneously study the nontrivial $\R^3$ bundle over $S^2$, whose associated sphere bundle
has total space $CP^2\#\overline{CP^2}$.

\begin{prop}\label{P:con} Every nonnegatively curved connection metric on an $\R^3$-bundle over $S^2$ is isometric
to a Riemannian quotient of the form:
$$M=((S^3,g_0)\times(\R^3,g_f))/S^1,$$
where $g_f$ is an $S^1$-invariant metric on $\R^3$ and $g_0$ is a connection metric on the principal
bundle $S^1\hookrightarrow S^3\ra S^2$.
\end{prop}

Notice that an integer $k$ determines the relative speeds at which $S^1$ acts on $S^3$ and $\R^3$, and the bundle is
trivial if and only if $k$ is even.  The proposition implies that the holonomy group of the normal bundle of the soul
(the ``normal holonomy group'') is either trivial or isomorphic to $S^1$.  The proposition is actually
a corollary of a similar fact for $S^2$-bundles:
\begin{prop}\label{P:s2} For an $S^2$-bundle $\pi:M\ra B=S^2$, suppose that $M$ and $B$ have nonnegatively curved
metrics so that $\pi$ is a Riemannian submersion with totally geodesic fibers.
Then $\pi:M\ra B$ can be metrically re-described as
$$M=((S^3,g_0)\times(S^2,g_1))/S^1\stackrel{\pi}{\ra}(S^3,g_0)/S^1,$$
where $g_0$ is a connection metric on the principal bundle $S^1\hookrightarrow S^3\ra S^2$, and $g_1$ is an $S^1$-invariant metric on
$S^2$.
\end{prop}

In the both propositions, $g_0$ is a connection metric on $S^3$,
which means that (1) the principal $S^1$ action is by isometries, and (2) all it's orbits have the same length.
If instead we assume only (1), then $g_0$ is called a ``warped connection metric'', and
the resulting metric on the $\R^3$ (or $S^2$) bundle over $S^2$ is a non-connection metric with holonomy group $S^1$.
However, not all nonnegatively curved metrics with holonomy group $S^1$ are of this type, as the next example
shows.

\begin{example}\label{E}  In Proposition~\ref{P:con}, if $g_f$ is a product metric, i.e., $(\R^3,g_f)=(\R^2,g_1)\times\R$,
then
$$M=((S^3,g_0)\times(\R^3,g_f))/S^1 = (((S^3,g_0)\times(\R^2,g_1))/S^1)\times\R = (E,g_E)\times\R,$$
where $(E,g_E)$ is the total space of an $\R^2$-bundle over $S^2$ with a connection metric.  Walschap
showed in~\cite[Theorem 2.1]{CodimTwo} that when $E$ is nontrivial ($k\neq 0$), such a connection metric $g_E$ can be non-rigid.
The most natural one,
coming from the round $g_0$ and the flat $g_1$, can be altered fairly arbitrarily away from the soul
without losing nonnegative curvature.  If $g_{E'}$ is such an alteration, then $(\R^2,g_{E'})\times\R$
has nonnegative curvature and normal holonomy group $S^1$.  Considering small metric spheres about a soul
of $(\R^2,g_{E'})\times\R$
produces a large family of nonnegatively curved metrics on $S^2\times S^2$ ($k$ even) and on $CP^2\#\overline{CP^2}$ ($k$ odd).
\end{example}

Although the metrics is the previous example are flexible, they are also rigid in the following sense:

\begin{prop}\label{P:quasi1}
For any complete nonnegatively curved metric on an $\R^3$-bundle over $S^2$ with normal holonomy group $S^1$,
the inequality of Proposition~\ref{P:condition} is not quasi-strict.
\end{prop}

The remaining case is when the normal holonomy group is transitive, and we begin by showing that examples
of this phenomenon exist:

\begin{prop}\label{P:ex}
Consider the following nonnegatively curved metric on $S^2\times\R^3$:
$$M=((S^2,\text{round})\times(\R^3,g_f)\times(\SO(3),g_B))/SO(3),$$
where $g_f$ is a nonnegatively curved $SO(3)$-invariant metric on $\R^3$ and $g_B$ is a bi-invariant metric on $\SO(3)$.
Then
\begin{enumerate}
\item The soul of $M$ is round.
\item The normal holonomy group is $\SO(3)$.
\item Each Sharafutdinov fiber is isometric to $\R^3$ with a fixed $S^1$-invariant metric.
\end{enumerate}
\end{prop}

If $g_B$ is replace by a right-invariant nonnegatively curved metric in this construction,
then $M$ is still nonnegatively curved.  These are the only known examples with transitive holonomy.  In
particular, for the non-trivial $\R^3$-bundle over $S^2$, it is not known whether there is a nonnegatively
curved metric with transitive normal holonomy group.
Although it is difficult to describe the general right-invariant case as explicitly as we describe the bi-invariant
case, it is at least straightforward to check that the metric spheres about the soul are not quasi-positively curved,
so the inequality of Proposition~\ref{P:condition} is not quasi-strict.

Since there are no other known examples with with transitive holonomy, one might ask whether an arbitrary metric
with transitive holonomy looks like with one of these examples, at least at the soul.  In section~\ref{S:s1}, we
show some special cases in which the inequality of Proposition~\ref{P:condition} rigidly determines what the metric may look
like at the soul.  For example, if we assume that the  Sharafutdinov fibers are all $S^1$-invariant,
then the inequality is not quasi-strict.  If we additionally assume that the connection
in the normal bundle of the soul has a special form, then the the soul must be round.
These added assumptions are motivated by the geometry of the known examples, and our rigidity statements
provide evidence that the family of nonnegatively curved metrics on $S^2\times\R^3$ with transitive holonomy is small.

The author is pleased to thank Detlef Gromoll for helpful conversations about these results.
\section{Connection metrics}
In this section, we prove Proposition~\ref{P:con}, which classifies the nonnegatively curved connection metrics
on an $\R^3$-bundle over $S^2$.

\begin{lem}\label{L:S1}For a connection metric of nonnegative curvature on an $\R^3$-bundle over $S^2$, the normal
holonomy group is trivial or is isomorphic to $S^1$.
\end{lem}
\begin{proof}
Let $M$ be the space $S^2\times\R^3$ together with a nonnegatively curved connection metric.  Let $\Sigma$
be a soul of $M$.  Let $p\in\Sigma$, and let $X,Y$ be an orthonormal basis of $T_p\Sigma$.
Let $\RN$ denote the curvature tensor of the connection in the normal bundle of $\Sigma$.

The map $U\ra\RN(X,Y)U$ is a skew-symmetric endomorphism of the fiber $\nbp=\R^3$, so there must be a non-zero
vector $W$ in it's kernel.  Choose any $V\in\nbp$ such that $\{V,W\}$ is orthonormal.
Let $\alpha:[0,a]\ra\Sigma$ be a piecewise-geodesic loop at $p$.  Let
$X_t,Y_t,W_t,V_t$ denote the parallel transports of $X,Y,W,V$ along $\alpha$.  From Proposition~\ref{P:condition},
$$\lb (D_{\alpha'(t)}\RN)(X_t,Y_t)W_t,V_t\rb^2\leq|\RN(W_t,V_t)X_t|^2\cdot R(X_t,Y_t,Y_t,X_t).$$
In other words, if we let $f(t)=|\RN(W_t,V_t)X_t| = \lb\RN(X_t,Y_t)W_t,V_t\rb$,
then
$$f'(t)^2\leq f(t)^2\cdot R(X_t,Y_t,Y_t,X_t).$$
Since $W$ is chosen so that $f(0)=0$, it follows that $f(t)=0$ for all $t$.  So either $W$ is fixed by the action
of the normal holonomy group on $\nbp$, or $\RN=0$, which implies that the normal holonomy group is trivial.
\end{proof}

The following would be a natural generalization of Lemma~\ref{L:S1}
\begin{conj} For a connection metric of nonnegative curvature on an $\R^k$-bundle over a manifold $M^n$,
the normal holonomy group is isomorphic to a subgroup of $\SO(n)$.
\end{conj}

Next we prove a result for $S^2$ bundles over $S^2$ analogous to Lemma~\ref{L:S1}.
\begin{lem} For an $S^2$-bundle $\pi:M\ra B=S^2$, suppose that $M$ and $B$ have nonnegatively curved
metrics so that $\pi$ is a Riemannian submersion with totally geodesic fibers.  Then the holonomy group of $\pi$
is trivial or is isomorphic to $S^1$.
\end{lem}
\begin{proof}
Let $\Hor,\V$ denote the horizontal and vertical distributions of $\pi$, and let $A$ denote the $A$-tensor of $\pi$.
Since $\Hor$ is 2-dimensional, $\U_p=A(X,Y)$ does not depend on the choice of oriented orthonormal basis $\{X,Y\}$
of $\Hor_p$.  Thus $\U$ is a well-defined global vertical vector field on $M$.

Let $p\in M$ be a point where $\U_p=0$.  There must be such a point on each $\pi$-fiber, since $S^2$ does not admit a
nowhere vanishing vector field.  Let $X,Y$ be an oriented orthonormal basis for $\Hor_p$,
and let $V\in\V_p$ be arbitrary.  Using O'Neill's formula for Riemannian submersions,
\begin{eqnarray*}
R(X,Y+V,Y+V,X) & = & R(X,Y,Y,X) + 2R(X,Y,V,X) + R(X,V,V,X)\\
               & = & R(X,Y,Y,X) + 2\lb (D_XA)_XY,V\rb + 0\\
               & = & R(X,Y,Y,X) + 2\lb \nabla_X\U,V\rb \geq 0
\end{eqnarray*}
Since the above is true for any $V$ (of any length), it follows that $(\nabla_X\U )^{\V}=0$.  It's also
easy to see that $(\nabla_X\U)^{\Hor}=0$ since for any $Z\in\Hor_p$,
$\lb\nabla_X\U,Z\rb = \lb A_XU,Z\rb = 0$.

Since the fibers are totally geodesic, $\U$ restricted to any fiber is a Killing field.  So, assuming $\U$ does not
vanish on the entire fiber $F_p$ through $p$, $p$ is an isolated zero of $\U$ on $F_p$.
In a neighborhood of $p$ in $M$, the set $S$ of zeros of $\U$ is a 2-dimensional smooth submanifold.  But,
$$T_pS \subset \{X\in T_pM|\nabla_X\U=0\} = \Hor_p.$$
So $T_pS=\Hor_p$, and $S$ is an integral submanifold of $\Hor$ near $p$.  This implies that $p$ is a fixed point of the holonomy
group.
\end{proof}

\begin{proof}[Proof of Proposition~\ref{P:s2}]
Choose a fiber of $\pi:M\ra B=S^2$.  Define $\U$ as in the previous proof.
$\U$ restricted to this fiber is a Killing field.  Let $O$ denote an orbit of maximal
length (the ``equator'' of the fiber).  Let $N$ denote the orbit of $O$ under the holonomy group, which
is a smooth 3-dimensional submanifold of $M$ (the union of the equators of all of the
fibers).  It is straightforward to verify that $N$ is totally geodesic.

The induced metric on the circle bundle $S^1\hookrightarrow N\stackrel{\pi|_N}{\ra} B$ is a connection metric, which we
claim can be metrically re-described as:
$$N=((S^3,g_0)\times S^1(r))/S^1\stackrel{\pi|_N}{\ra} B=(S^3,g_0)/S^1.$$
More precisely, the relative speed, $k$, at which $S^1$ acts on $S^3$ and $S^1$ can be first chosen to give
the topologically correct bundle.  Next, choosing a connection metric $g_0$ on the bundle $S^1\hookrightarrow S^3\ra S^2$ means choosing
(1) a metric on $S^2$, (2) a principal connection in the bundle, (3) the fiber-length, $l$.  We make the first choice
such that $(S^3,g_0)/S^1$ is isometric to $B$.  We make the second choice to induce the correct horizontal distribution
in the circle bundle.  Finally, we can choose $l$ and $r$ together so that the fiber-length is correct.

If $g_1$ denotes any $S^1$-invariant metric on $S^2$ whose equator (meaning the maximal-length
orbit of the $S^1$-action) has circumference $2\pi r$, then the Riemannian submersion
$$\tilde{M}=((S^3,g_0)\times(S^2,g_1))/S^1\stackrel{\tilde{\pi}}{\ra}B=(S^3,g_0)/S^1$$
has totally geodesic fibers and holonomy group $S^1$.  Further, if $\tilde{N}$ denotes
the union of the equators of all of the fibers of $\tilde{\pi}$, then $\tilde{N}\stackrel{\tilde{\pi}|_{\tilde{N}}}{\ra} B$ is metrically equal to
$N\stackrel{\pi|_N}{\ra} B$.  Finally, $g_1$ can easily be chosen so that the fiber metric of $\tilde{\pi}$ agrees with
the fiber metric on $\pi$.  It's then straightforward to see that $\tilde{M}\stackrel{\tilde{\pi}}{\ra} B$
is metrically equal to $M\stackrel{\pi}{\ra} B$.
\end{proof}
\begin{proof}[Proof of Proposition~\ref{P:con}]
For a nonnegatively curved connection metric on an $\R^3$-bundle, $M$, over $S^2$, the distance sphere, $S_r$,
of any radius $r$ about a soul inherits a metric which, by Proposition~\ref{P:s2}, can be described as:
$$S_r=((S^3,g_0)\times(S^2,g_1))/S^1.$$
As we vary $r$, the soul metric $(S^3,g_0)/S^1$ does not change, and neither does the
horizontal distribution of $(S^3,g_0)\ra(S^3,g_0)/S^1$, since horizontal displacement of points in the sphere
bundles is controlled by parallel transport of vectors in the normal bundle of the soul.  So we can choose
$g_0$ and $k$ independent of $r$.  We write $g_1(r)$ to show the dependence of $g_1$ on $r$.  Now define
an $S^1$-invariant metric $g_f$ on $\R^3$ such that the distance sphere of radius $r$ about the origin
has metric $g_1(r)$.  It follows that $M$ is isometric to $((S^3,g_0)\times(\R^3,g_f))/S^1.$
\end{proof}
\section{Non-connection metrics with normal holonomy group $S^1$}
We saw in example~\ref{E} that there is no simple group-theoretic classification of nonnegatively
curved metrics on $\R^3$-bundles over $S^2$ with normal holonomy group $S^1$.
In this section, we prove Proposition~\ref{P:quasi1}, which says that even though
such metrics are in one sense perturbable, there is always partial rigidity at the soul.
\begin{proof}[Proof of Proposition~\ref{P:quasi1}]
Let $M$ denote the total space of an $\R^3$-bundle over $S^2$ together with a complete metric of nonnegative
curvature.  Let $\Sigma$ denote a soul of $M$.
Assume that the normal holonomy group is isomorphic to $S^1$, so that there exists a parallel section, $W$,
of the normal bundle, $\nb$, of $\Sigma$.

For any $p\in\Sigma$, orthonormal vectors $X,Y\in T_p\Sigma$ and any unit-length vector $V\in\nbp$ orthogonal to $W$,
we have:
$$\lb (D_XR)(X,Y)W,V\rb = |R(W,V)X| = 0.$$
So the inequality of Proposition~\ref{P:condition} implies that $(D_{X,X}R)(W,V,V,W)\geq 0$.  We claim that in fact,
$(D_{X,X}R)(W,V,V,W)=0$, which suffices to demonstrate that the inequality is not quasi-strict.
To prove this claim, we must show that the curvature of any 2-plane normal to the soul
containing $W$ is the same as any other.

To see this, assume $p$ and $V$ were chosen so that $R(W,V,V,W)$ is maximal among all possible choices of $p$ and $V$.
Let $\alpha(t)$ denote any piecewise geodesic in $\Sigma$ with $\alpha(0)=p$,
and let $X_t,Y_t,W_t,V_t$ be the parallel transports of $X,Y,W,V$ along $\alpha$.
Let $f(t)=R(W_t,V_t,V_t,W_t)$.  The inequality says that $f''(t)\geq 0$.  By the maximality assumption,
$f(t)$ must be a constant function.  But since the holonomy group is $S^1$,
the 2-planed at $p$ spanned by $W$ and $V$ can be parallel transported along piecewise geodesics
arbitrarily close to any other vertical 2-plane containing
$W$ at any other point of the soul. Therefore, both the left and right sides of the inequality
become zero when one of the vertical vectors equals $W$.
\end{proof}

\section{An example with transitive normal holonomy group}\label{S:example}
In this sections, we prove Proposition~\ref{P:ex}, which describes the geometry of the Riemannian manifold
$$M=((S^2,\text{round})\times(\R^3,g_f)\times(\SO(3),g_B))/SO(3),$$
To simplify the discussion, we take the unit-round metric on $S^2$,
we take $g_f$ to be the flat metric on $\R^3$, and we assume the bi-invariant metric $g_B$ is scaled so that
a unit-length vector in $\so(3)$ corresponds to a Killing field on $S^2(1)$ with maximal norm one.  Let
$\Soul$ denote the soul of $M$.
We begin by explicitly describing parallel transport in $\nb$, proving in particular:
\begin{lem}\label{L:e1}
The normal holonomy group of $M$ is $SO(3)$.
\end{lem}
\begin{proof}
Let $g$ denote the product of the unit-round
metric on $S^2$ with the flat metric on $\R^3$.  Let $\tilde{g}$ denote the quotient metric on
$S^2\times\R^3$ obtained from the above description of $M$.
According to~\cite{C}, $\tilde{g}$ is obtained from $g$ by ``rescaling along the $SO(3)$
action.''  We next describe more precisely how the metric $\tilde{g}$ on $T_{(p,V)}(S^2\times\R^3)$
is obtained from the metric $g$.

Let $(p,V)\in S^2(1)\times\R^3$ with $|V|=1$.  Using the natural inclusion $S^2(1)\subset\R^3$,
define $\theta$ as the angle between $p$ and $V$.  Assume that $V\neq \pm p$. Denote:
$$W=\frac{V-\lb V,p\rb\cdot p}{|V-\lb V,p\rb\cdot p|}\text{\,\,\, and \,\,\,} U=W\times p,$$
where ``$\times$'' denotes the vector cross product in $\R^3$.  Notice that $\{W,p,U\}$ is an
oriented orthonormal basis of $\R^3$, and $V=\cos\theta\cdot p + \sin\theta\cdot W$.
Next, consider the following orthonormal basis $\{X,A,Y,B,\hat{r}\}$ of $T_{(p,V)}(S^2\times\R^3)$:
$$X=(-W,0), A=(U,0), Y=(0,\sin\theta\cdot p - \cos\theta\cdot W), B=(0,U), \hat{r}=(0,V).$$


Choose an orthonormal basis $\{E_1,E_2,E_3\}$ of $\text{so}(3)$ corresponding to unit-speed
right-handed rotations of $S^2$ about the vectors $\{U,p,W\}$ respectively.
The values, $T_1,T_2,T_3\in T_{(p,V)}(S^2\times\R^3)$, of the corresponding Killing fields at
$(p,V)$ are:
$$T_1=X+Y,\text{\,\,\,\,\,\,}T_2=\sin\theta\cdot B,\text{\,\,\,\,\,\,}T_3=A+\cos\theta\cdot B.$$
Define $$K=\left( k_{ij} \right) = \left( \lb T_i,T_j\rb_{g}\right) = \left( \begin{array}{ccc}
2 & 0 & 0 \\ 0 & \sin^2\theta & \sin\theta\cdot\cos\theta \\ 0 & \sin\theta\cdot\cos\theta & 1 +
\cos^2\theta \end{array}\right).$$
The metric $\tilde{g}$ agrees with $g$ on the orthogonal compliment of $\text{span}\{T_i\}$.  On
$\text{span}\{T_i\}$ we have:
$$\tilde{K}=\left( \tilde{k}_{ij} \right) = \left( \lb T_i,T_j\rb_{\tilde{g}}\right) = K\cdot(I+K)^{-1} =
\left( \begin{array}{ccc}
\frac{2}{3} & 0 & 0 \\ 0 & 2\frac{\cos^2\theta - 1}{\cos^2\theta - 4} &
-\frac{\sin\theta\cdot\cos\theta}{\cos^2\theta - 4} \\ 0 & -\frac{\sin\theta\cdot\cos\theta}{\cos^2\theta - 4} &
\frac{-2}{\cos^2\theta - 4}\end{array}\right).$$

It is useful to define $T_4=X-Y$, so that $\{T_1,T_2,T_3,T_4,\hat{r}\}$ becomes a basis for
$T_{(p,V)}(S^2\times\R^3)$.  The change of basis matrix which translates from
$\{X,A,Y,B,\hat{r}\}$ to $\{T_1,T_2,T_3,T_4,\hat{r}\}$ is $N = \left(
\begin{array}{ccccc} \half & 0 & \half & 0  & 0 \\ 0 & -\frac{\cos\theta}{\sin\theta} & 0 &
\frac{1}{\sin\theta} & 0  \\ 0 & 1 & 0  & 0 & 0 \\  \half & 0 & -\half & 0 & 0
\\ 0 & 0 & 0 & 0 & 1\end{array}\right).$
So if $\tilde{K}'=\left(\begin{array}{ccc}\tilde{K} & 0 & 0 \\ 0 &
2 & 0 \\ 0 & 0 & 1 \end{array}\right)$, then two vector $\mathcal{X},\mathcal{Y}\in
T_{(p,V)}(S^2\times\R^3)$ written in the basis $\{X,A,Y,B,\hat{r}\}$ have inner product:
\begin{equation}\label{metric}
\lb \mathcal{X},\mathcal{Y}\rb_{\tilde{g}} = (N\cdot\mathcal{X})^T\cdot\tilde{K}'\cdot (N\cdot\mathcal{Y}).
\end{equation}

Since the metric in the $\hat{r}$ direction is unchanged when passing from $g$ to $\tilde{g}$,
the Sharafutdinov map $\pi:M\ra\Soul$ is also unchanged when passing from $g$ to $\tilde{g}$.  So the
vertical space of $\pi$ is unchanged:
$$\V_{(p,V)} = \text{span}\{Y,B,\hat{r}\} = T_V\R^3\subset T_{(p,V)}(S^2\times\R^3).$$
It is now easy to show that the horizontal space of $\pi$ is
\begin{equation}\label{hor}
\Hor_{(p,V)} = \text{span}\{X+\half\cdot Y,A+\half\cos\theta\cdot B\}.
\end{equation}
To verify this, use equation~\ref{metric} to show that both of these vectors are
$\tilde{g}-$perpendicular to $Y$, $B$, and $\hat{r}$.

Equation~\ref{hor} allows a simple description of $\tilde{g}$-parallel transportation in $\nb$.  Namely,
suppose that $\gamma$ is a geodesic segment in $\Soul$ with $\gamma(0)=p$.  Using the $g$-parallel
identification $\R^3=\nbp=\nu_{\gamma(t)}(\Soul)$, we can think of the parallel transport of $V\in\nbp$ along
$\gamma$ as a path, $V(t)$, in $\R^3$.  If $\mathcal{U}\in\Soul$ is perpendicular to the plane in which $\gamma$ lies,
then $V(t)$ is the result of rotating $V$ an angle of $t/2$ about the axis determined by $\mathcal{U}$.  In other words,
parallel transport rotates $\R^3$ along with $\gamma$, but at half the speed.  The Lemma follows easily.
\end{proof}

\begin{lem}\label{L:e2}
The soul is a round sphere with constant curvature $2$.
\end{lem}
\begin{proof}
As before, let $T_1,T_2,T_3\in T_{(p,0)}(S^2\times\R^3)$ denote the values at $(p,0)$ of the Killing fields associated with
an orthonormal basis $E_1,E_2,E_3$ of $\text{so}(3)$.  If the $E_i$'s are chosen as before, then
$T_2=0$ and $\{T_1,T_3\}$ form a $g$-orthonormal basis of $T_pS^2\subset T_{(p,0)}(S^2\times\R^3)$.
So $\lb T_1,T_3\rb_{\tilde{g}}=0$ and for $i=1,3$, $|T_i|_{\tilde{g}}=|T_i|_g/\sqrt{1+|T_i|_g}=1/\sqrt{2}$.  In other words,
the metric on the soul is the unit-round metric with the norms of all vectors rescaled by a factor of $1/\sqrt{2}$.
This metric has constant curvature $2$.
\end{proof}
 
\begin{lem}\label{L:e3}
The intrinsic metric on a fiber of the Sharafutdinov map $\pi:M\ra\Soul$ is $S^1$-invariant.  Any two such
fibers are isometric.
\end{lem}

\begin{proof}
Using notation from the proof of Lemma~\ref{L:e1}, a $g$-orthonormal basis for the tangent space to the
fiber $F_p=\pi^{-1}(p)$ at $(p,V)$ is $\{Y,B,\hat{r}\}$.  By equation~\ref{metric},
$$ \lb Y,Y\rb_{\tilde{g}} = 2/3, \qquad \lb Y,B\rb_{\tilde{g}} = 0, \qquad \lb B,B\rb_{\tilde{g}} = \frac{2}{4-\cos^2\theta}.$$
Therefore, the metric on the distance sphere of radius $1$ about $p$ in
$F_p$ is $S^1$-invariant and independent of $p$.  But nothing essential changes in the proof of Lemma~\ref{L:e1}
if you take $V$ to have norm other
than 1, so distance spheres of other radii about $p$ in $F_p$ are also $S^1$-invariant, so
$F_p$ is $S^1$-invariant.  The fixed vector of the $S^1$-symmetry is $p$.
\end{proof}

The $S^1$-symmetry of the fibers implies that the curvature of a vertical 2-plane $\sigma$ at a point $p\in\Soul$
depends only the angle $\theta$ that $\sigma$ makes with $p$.  We leave it to the reader to verify
that planes containing $p$ have curvature $3$, while the plane orthogonal to $p$ has curvature $3/2$, so
\begin{equation}\label{vert} k(\sigma)=3\sin^2\theta + \frac{3}{2}\cos^2\theta.\end{equation}

Next we describe the curvature tensor $\RN$ of the connection $\nabla$ in the normal bundle of the soul.
\begin{lem}\label{L:e4} Let $p\in\Soul$, let $\{X,Y\}$ be an oriented orthonormal basis of $T_p\Soul$, and
let $V\in\nbp$.  Letting ``$\times$'' denote the vector cross product in $\R^3$, and using the
natural identification $\Soul=S^2\subset\R^3=\nbp$, we have:
$$\RN(X,Y)V=\frac{3}{2} V\times p.$$
In particular, $|\RN|=\frac{3}{2}$ and $\RN(X,Y)p=0$.
\end{lem}
Although we could prove Lemma~\ref{L:e4} directly from equation~\ref{hor}, we find it simpler
and more illuminating to prove a more general formula in the next section.


\section{A family of connections in a trivial $\R^3$-bundle}\label{S:family}
In this section, we study a natural family of connections in the trivial $\R^3$-bundle over a manifold $B^n$.

Let $f:B^n\rightarrow S^2$ denote any smooth function.  Let $\lambda$ be any real number.
There is a connection, $\nabla$, in the trivial vector bundle $B^n\times\R^3$ naturally
associated with $\{f,\lambda\}$ as follows.  Let $\Phi:TS^2\rightarrow so(3)$ denote the canonical map,
which can be described algebraically as
$$\Phi(p,X)(W) = (p\times X)\times W = \lb p,W\rb X - \lb X,W\rb p.$$
Then $\mathcal{D}=\Phi\circ df:TB^n\rightarrow so(3)$ is a connection difference form; that is,
$\mathcal{D}$ is a $so(3)$-valued one-form on $B^n$.  So if $\ON$ denotes the flat
connection on $B^n\times\R^3$, then $$\nabla = \ON + \lambda\mathcal{D}$$ is a connection on $B^n\times \R^3$.

\begin{lem}
The connection in the normal bundle of the soul of $$((S^2,\text{round})\times(\R^3,g_f)\times(\SO(3),g_B))/SO(3)$$
is determined as above by the identity map $f:S^2\ra S^2$, with $\lambda=-1/2$.
\end{lem}
\begin{proof}
This follows immediately from equation~\ref{hor}.
\end{proof}

Thinking of $f$ as a unit-length section of the bundle, the next lemma says that the curvature tensor, $\RN$,
of $\nabla$ vanishes along $f$, and the norm of $\RN$ depends on the area-distortion of $f$.
Notice that Lemma~\ref{L:e4} is a corollary of Lemma~\ref{L:curvL} and Lemma~\ref{L:e2}.

\begin{lem}\label{L:curvL}
Let $p\in B^n$, $X,Y\in T_p B^n$ and $W\in\R^3$.  Define $\BX=df_p X$ and $\BY=df_p Y$.
Then
$$\RN(X,Y)W = \lambda(\lambda+2)\left(\lb\BX,W\rb \BY - \lb\BY,W\rb\BX\right).$$
In particular, $\RN(X,Y)f(p)=0$ and $|\RN(X,Y)|= |\lambda^2+2\lambda|\cdot|\BX\wedge\BY|$.
\end{lem}
\begin{proof}
Extend $\{X,Y\}$ to local vector fields on $B^n$ which commute at $p$.  Extend $W$ to a
$\overline{\nabla}$-parallel section of the bundle.  Let $\bp=f(p)$.  Then
\begin{eqnarray*}
\RN(X,Y)W & = & \nabla_X\nabla_Y W - \nabla_Y\nabla_X W \\
          & = & \lambda(\nabla_X\D_Y W - \nabla_Y\D_XW) \\
          & = & \lambda(\ON_X\D_YW-\ON_Y\D_XW) + \lambda^2(\D_X\D_YW - \D_Y\D_X W).
\end{eqnarray*}
We have at $p$ that:
$$\D_Y W = (\bp\times\BY)\times W = \lb \bp,W\rb\BY - \lb\BY,W\rb\bp,$$
so
$$\D_X\D_Y W = \lb \bp,\D_YW\rb\BX - \lb\BX,\D_YW\rb \bp = -\lb\BY,W\rb\BX - \lb \bp,W\rb\lb\BX,\BY\rb p,$$
and therefore,
$$\D_X\D_Y W-\D_Y\D_X W = \lb\BX,W\rb\BY - \lb\BY,W\rb\BX.$$

Next let $p(t)$ denote a path in $B^n$ with $p(0)=p$ and $p'(0)=X$.
It is convenient to let $\bar{p}(t)=f(p(t))$ and $\BY(t)=df_{p(t)} (Y(p(t)))$.  then,
$$(\D_YW)_{p(t)} = (\bar{p}(t)\times\BY(t))\times W = \lb \bar{p}(t),W\rb\BY(t) - \lb \BY(t),W\rb\bar{p}(t)$$
So at $p$,
\begin{eqnarray*}
\ON_X \D_YW & = & \frac{D}{dt}\bigg|_{t=0}
    \left( \lb \bar{p}(t),W\rb\BY(t) - \lb \BY(t),W\rb\bar{p}(t)\right)\\
            & = & \lb\BX,W\rb\BY +\lb\bp,W\rb\BY'(0)-\lb\BY'(0),W\rb\bp - \lb\BY,W\rb\BX.
\end{eqnarray*}
Therefore,
\begin{eqnarray*}
\ON_X \D_YW - \ON_Y \D_X W & = & 2\lb\BX,W\rb\BY - 2\lb\BY,W\rb\BX \\
                           &   & +\lb\BX'(0)-\BY'(0),W\rb\bp - \lb\bp,W\rb(\BX'(0)-\BY'(0)) \\
                           & = & 2\lb\BX,W\rb\BY - 2\lb\BY,W\rb\BX
\end{eqnarray*}

The last equality is justified because $[X,Y](p)=0$ implies that $[\BX,\BY](\bp)=0$, which in turn implies
that $\BX'(0)-\BY'(0)$ is parallel to $\bp$.  This completes the proof.
\end{proof}
 
We end this section by mentioning a few interesting properties of the connection associated with $\{f,\lambda\}$.
If $\lambda=0$, then $\nabla$ is clearly flat.  By Lemma~\ref{L:curvL}, if $\lambda=-2$
then $\nabla$ is flat, which is perhaps less obvious.  Regarding $f$ as a unit-length
section of the bundle, the covariant derivative of $f$ in the direction $X\in T_p B^n$ is
\begin{equation}\label{DXW}
\nabla_X f = \ON_X f + \lambda\D_X f(p) = \BX + \lambda (f(p)\times\BX)\times f(p) = \BX + \lambda\BX = (1+\lambda)\BX,
\end{equation}
where $\BX=df_pX$.  So if $\lambda=-1$, then $f$ is a parallel section, and hence the holonomy group of $\nabla$ is
isomorphic to $S^1$.
Thus, $\lambda$ plays a significant role along with $f$ in determining the qualitative geometric properties of $\nabla$.


\section{Rigidity with $S^1$-invariant Sharafutdinov fibers}\label{S:s1}

In this section, we show that for a class of metrics more general than connection metrics,
the inequality of Proposition~\ref{P:condition} forces some rigidity for the metric at the soul.
As before, $M$ will denote the space $S^2\times\R^3$ together with a metric of nonnegative curvature,
$\Soul$ will denote a soul of $M$, and $\nabla$ will denote the connection in the normal bundle $\nb$.
In the remainder of this section, we make the following assumptions:
\begin{enumerate}
\item $M$ has a curvature nullity section, i.e., a global unit-length section $W$ of $\nb$ such that
$\RN(\cdot,\cdot)W(p)=0$ for all $p\in\Soul$.
\item The sectional curvature of a 2-plane $\sigma\subset\nbp$ depends only on the angle that $\sigma$
forms with $W(p)$.
\end{enumerate}

The first assumption is true in all known examples.  To understand its content, define $F:\Soul\ra\R$ as the norm of $\RN$.
Let $p\in\Soul$ be a point at which $F(p)\neq 0$.
Let $X,Y\in T_p\Soul$ be an oriented orthonormal basis.  Since $\RN(X,Y):\nbp\ra\nbp$ is a skew-symmetric endomorphism
of a $3$-dimensional vector space, there is a unit-length vector $W(p)\in\nbp$ such that $\RN(X,Y)W(p)=0$.
If $F>0$ on $\Soul$, then a global unit-length section of $\nb$, $p\mapsto W(p)$, can be constructed.
In this case, there is a complimentary bundle $$\R^2\hookrightarrow W^{\perp}=\{U\in\nb:\lb U,W\rb =0\}\ra\Soul.$$
Even in known examples where $F$ is not strictly positive, it is always possible to find a curvature nullity
section $W$, and therefore to define $W^{\perp}$.  In general, $W$ is not a parallel section.

For example, if $\nabla$ is determined by $\{f:S^2\ra S^2,\lambda\}$ as in section~\ref{S:family}, then
$f$ is a curvature nullity section, and the
isomorphism class of the complimentary bundle depends on the mapping degree of $f$.  If $f$ is a diffeomorphism,
then the complimentary bundle is isomorphic to $TS^2$.
As a second example, if $M$ has normal holonomy group $S^1$, then there is a parallel curvature nullity
section.  In the connection metrics of Proposition~\ref{P:con}, the isomorphism class of the complimentary
bundle depends on the even integer $k$.

The second assumption says that each fiber looks at the soul
as if it admits an isometric $S^1$ actions with fixed direction $W$.  This is true in the example
of Section~\ref{S:example}.  If $U,V\in\nbp$ are orthonormal
vectors orthogonal to $W=W(p)$, the second assumption implies
$R(U,V)W=R(U,W)V=0$. It also implies that $g_0(p)=R(W,U,U,W)$ and $g_1(p)=R(U,V,V,U)$ describe
well-defined functions $g_0,g_1:\Soul\ra\R$.  Compare with equation~\ref{vert}, where
$g_0$ and $g_1$ are constant.  Proposition~\ref{P:con} implies the following
relationships between these functions:

\begin{prop}\label{three}
Let $p\in\Soul$ and let $X\in T_p\Soul$ be unit-length.  Define $a(X)=|\nabla_XW|$, and let $k_{\Soul}$
denote the Gauss curvature of $\Soul$.  Then,
\begin{gather*}
(XF)^2 \leq k_{\Soul}\left(F^2 + \frac{2}{3}\text{hess}_{g_1}(X,X) -
\frac{4a(X)^2}{3}(g_1-g_0)\right)\label{g1} \\ a(X)^2\cdot F^2 \leq
\frac{2}{3}k_{\Soul}\left(\text{hess}_{g_0}(X,X) + 2a(X)^2(g_1-g_0)\right)\label{g2} \\ 0 \leq
\frac{2}{3}k_{\Soul}\left(\text{hess}_{g_0}(X,X)\right)\label{g3} \end{gather*}
\end{prop}

\begin{proof}
Let $p\in\Soul$ and $X\in T_p\Soul$ with $|X|=1$.
Let $V=(\nabla_XW)/a(X)$ if $a(X)\neq 0$; otherwise
let $V\in\nbp$ be an arbitrary unit-length vector orthogonal to $W$.  Let $U\in\nbp$ be such that
$\{U,V,W\}$ is an orthonormal basis of $W$.  Choose $Y\in T_p\Soul$ so that
$\{X,Y\}$ forms an oriented orthonormal basis.  For $E_1,E_2\in\nbp$, Define $G(X,E_1,E_2)$ as
$$k_{\Soul}\cdot(|\RN(E_1,E_2)X|^2+(2/3)(D_{X,X}R)(E_1,E_2,E_2,E_1))-\lb (D_X\RN)(X,Y)E_1,E_2\rb^2,$$
which is the right side minus the left side of the inequality of Proposition~\ref{P:condition}.
The inequalities of the proposition come from: $G(X,U,V)\geq 0$, $G(X,U,W)\geq 0$ and $G(X,V,W)\geq 0$.

Notice that $|\RN(W,U)X|=|\RN(W,V)X|=0$, while $|\RN(U,V)X|=F$.  Next, extending $U$ and $V$ to be parallel
along the path in the direction of $X$, it's easy to see that:
\begin{eqnarray*}
\lb(D_X\RN)(X,Y)W,V\rb & = & -\lb\RN(X,Y)(\nabla_XW),V\rb = 0 \\
\lb(D_X\RN)(X,Y)W,U\rb & = & -\lb\RN(X,Y)(\nabla_XW),U\rb = -a(X)\cdot F
\end{eqnarray*}
It remains to compute $\lb(D_X\RN)(X,Y)U,V\rb$ and the three terms involving
$(D_{X,X}R)$. For this we must find nice extensions of the vectors $\{U,V,W\}$ in the
direction $X$.  Let $\alpha(t)$ be the geodesic in $\Soul$ with $\alpha(0)=p$ and $\alpha'(0)=X$.
Let $W(t)=W(\alpha(t))$.  For a first try, let $U(t)$ and $V(t)$ be arbitrary extensions of
$U$ and $V$ along $\alpha(t)$ so that $\{U(t),V(t),W(t)\}$ is orthonormal for all $t$.

By choice of $V=V(0)$, we have $W'(0)=a\cdot V(0)$, where $a=a(X)$.  Also,
$ V'(0) = -a\cdot W + b\cdot U$ and $U'(0)=-b\cdot V$
for some $b$ depending on the choice of extensions.
We can improve our choice of extensions to achieve $b=0$ by defining:
\begin{eqnarray*}
\tilde{U}(t) & = &  \cos(\phi(t))\cdot U(t) + \sin(\phi(t))\cdot V(t) \\
\tilde{V}(t) & = & -\sin(\phi(t))\cdot U(t) + \cos(\phi(t))\cdot V(t),
\end{eqnarray*}
where $\phi(0)=0$, $\phi'(0)=b$, and
$\phi''(0)=\lb V''(0),U\rb = -\lb U''(0),V\rb$.  It is straightforward to compute that
$\tilde{U}'(0)=0$, $\tilde{V}'(0)=-aW$, $\tilde{U}''(0)\perp V$, and $\tilde{V}''(0)\perp U$.
To simplify notation, omit the tilde's.  Using these extensions,
$$\lb(D_X\RN)(X,Y)U,V\rb = XF - \lb\RN(X,Y)U'(0),V\rb - \lb\RN(X,Y)U,V'(0)\rb = XF.$$
The proposition will follow once we verify that:
\begin{eqnarray*}
(D_{X,X}R)(U,V,V,U) & = & \text{hess}_{g_1}(X,X) - 2a^2(g_1-g_0)\\
(D_{X,X}R)(W,U,U,W) & = & \text{hess}_{g_0}(X,X) + 2a^2(g_1-g_0)\\
(D_{X,X}R)(V,W,W,V) & = & \text{hess}_{g_0}(X,X)
\end{eqnarray*}
Since the arguments are similar, we verify the third equality, leaving the first two to the reader.
\begin{eqnarray*}
(D_{X,X}R)(V,W,W,V)
   & = & X((D_XR)(V,W,W,V)) - 2(D_XR)(V'(0),W,W,V) \\
   &   & -2(D_XR)(V,W'(0),W,V) \\
   & = & X((D_XR)(V,W,W,V)) - 0 - 0 \\
   & = & \text{hess}_{g_0}(X,X) - 2XR(V',W,W,V)-2XR(V,W',W,V),
\end{eqnarray*}
Where the last two terms simplify as follows:
\begin{eqnarray*}
XR(V',W,W,V) & = & (D_XR)(V'(0),W,W,V) + R(V''(0),W,W,V) \\
               &   &  + R(V'(0),W'(0),W,V) + R(V'(0),W,W'(0),V)\\
               &   & + R(V'(0),W,W,V'(0)) \\
               & = & 0 + \lb V''(0),V\rb g_0 -a^2R(W,V,W,V) + 0 + 0 \\
               & = & -a^2 g_0 + a^2 g_0 = 0 \\
XR(V,W',W,V) & = & (D_XR)(V,W'(0),W,V) + R(V'(0),W'(0),W,V) \\
               &   &  + R(V(0),W''(0),W,V) + R(V,W'(0),W'(0),V)\\
               &   & + R(V,W'(0),W,V'(0)) \\
               & = & 0 - a^2R(W,V,W,V)+ \lb W''(0),W\rb g_0  + 0 + 0 \\
               & = & a^2 g_0 - a^2 g_0 = 0
\end{eqnarray*}
\end{proof}
From the third inequality of Proposition~\ref{three} we immediately learn the following, with $G$ defined
as in the previous proof:
\begin{cor} $g_0$ is a constant function, and $G(X,W,V)=0$.  Hence, the inequality of Proposition~\ref{P:condition}
is not quasi-strict.
\end{cor}
In the example of section~\ref{S:example}, the first and second inequalities of Proposition~\ref{three}
are both strictly satisfied.  In fact, $\text{span}\{W,V\}$ is the unique plane for which $G(X,W,V)=0$.  We therefore
do not expect any rigidity to follow from the first and second inequalities.  However, we do get further
rigidity from the fact that $G(X,\cdot,\cdot)\geq 0$ for planes very near $\text{span}\{W,V\}$.
We use this idea to prove:
\begin{lem}\label{L:W} Let $\alpha(t)$ be a geodesic in $\Soul$.  Let $W(t)=W(\alpha(t))$ denote
the curvature nullity section restricted to $\alpha$.  Then
$W''(t)\in\text{span}\{W(t),W'(t)\}$ for all $t$.
\end{lem}
The lemma means that, even though $W$ is not parallel, there is a parallel 2-plane containing
$W$ along any geodesic; namely, the plane spanned by $W$ and $W'$.
\begin{proof}
Let $p=\alpha(0)$ and $X=\alpha'(0)$, which we can assume is unit-length.  Define $U,V\in\nbp$ as in the proof
of Proposition~\ref{three}, so that $V$ is parallel to $W'(0)$.
An arbitrary 2-plane $\sigma\subset\nbp$ is spanned by orthonormal vectors $E_1,E_2$ of the form:
\begin{eqnarray*}
E_1 & = & a_1W+b_1U+c_1V = (\cos\theta) W + 0U + (\sin\theta) V \\
E_2 & = & a_2W+b_2U+c_2V = (-\sin\phi\sin\theta)W + (\cos\phi)U + (\sin\phi\cos\theta)V
\end{eqnarray*}
Here, $E_2$ is an arbitrary unit-vector expressed in spherical coordinates, and $E_1$ is the unit-vector in the
$WV$-plane orthogonal to $E_2$.  Defining $G$ as in the previous proof,
\begin{eqnarray*}
\lefteqn{G(X,E_1,E_2)}\\ & = & (b_1c_2-b_2c_1)^2 G(X,U,V) + (a_1b_2-a_2b_1)^2 G(X,W,U) \\
           &   & -2(b_1c_2-c_1b_2)(a_1b_2-a_2b_1)\lb(D_X\RN)(X,Y)U,V\rb^2\lb(D_X\RN)(X,Y)W,U\rb^2 \\
           &   & -(4/3)b_1b_2(a_1-a_2)(c_2-c_1)(D_{X,X}R)(U,V,U,W)\\
           &   & -(4/3)c_1c_2(a_1-a_2)(b_2-b_1)(D_{X,X}R)(U,V,W,V)\\
           &   & -(4/3)(a_1^2c_2^2+a_2^2c_1^2-2a_1a_2c_1c_2)(D_{X,X}R)(W,V,V,W).
\end{eqnarray*}
Denote $G(\phi,\theta)=G(X,E_1,E_2)$.  Notice that $G(\pi/2,\theta)=G(X,W,V)=0$.  It's straightforward to
compute that:
$$\frac{dG}{d\phi}(\pi/2,\theta) = -(\sin\theta\cos^2\theta+\cos\theta-\cos^3\theta)\cdot (D_{X,X}R)(U,V,W,V)\geq 0.$$
Using an argument similar to the proof of Proposition~\ref{three}, we see that:
$$(D_{X,X}R)(U,V,W,V) = (g_1-g_0)\lb W''(0),U\rb.$$
The trigonometric expression changes sign as $\theta$ varies, so the only possibility is that
$(g_1-g_0)\lb W''(0),U\rb=0$.  But we can see from Proposition~\ref{three} that $g_1-g_0>0$ on $\Soul$,
so $\lb W''(0),U\rb = 0$,
which completes the proof.
\end{proof}

Lemma~\ref{L:W} has a strong corollary in the special case where the connection in
the normal bundle of the soul is of the type described in section~\ref{S:family}:
\begin{cor}
If $\nabla$ is induced by some diffeomorphism $f:S^2\ra S^2$ and some $\lambda\in\R$,
as described in section~\ref{S:family}, then $\Soul$ is round and $f$ is the identity function.
\end{cor}
\begin{proof}
Let $\alpha(t)$ be a geodesic in $\Soul$, and consider $W(t)=W(\alpha(t))=f(\alpha(t))$, which we could
think of as a section of $\nb$ along $\alpha$ or as a path in $S^2(1)$.
Using equation~\ref{DXW}, taking one more derivative, and comparing to Lemma~\ref{L:W}, we get that
$f(\alpha(t))$ has no geodesic curvature.  In other words, $f(\alpha(t))$ is a (possibly reparameterized) geodesic
on the round sphere $S^2(1)$.  Since $f:\Soul\ra S^2(1)$ maps geodesics to paths whose images are geodesics, it follows
that $\Soul$ must be a round sphere of some radius, and $f$ must be the identity map from $S^2$ to $S^2$.
\end{proof}

\bibliographystyle{amsplain}

\begin{thebibliography}{9}
\bibitem{C}
   J. Cheeger, \emph{Some examples of manifolds of nonnegative curvature}, J. Differential Geom.
   \textbf{8} (1972), 623-628.
\bibitem{GT} D. Gromoll and K. Tapp, \emph{Nonnegatively curved metrics on $S^2\times\R^2$}, Geometriae Dedicata,
   to appear.
\bibitem{LuisTestTube}
  L. Guijarro, \emph{Improving the metric in an open manifold with nonnegative curvature},
  Proc. Amer. Math. Soc. \textbf{126}, No. 5 (1998), 1541-1545.
\bibitem{GW}
  L. Guijarro and G. Walschap, \emph{The metric projection onto the soul},
  Tans. Amer. Math. Soc. \textbf{352} (2000), no. 1, 55-69.
\bibitem{SW}
   M. Strake and G. Walschap, \emph{Connection metrics of nonnegative curvature on vector bundles},
   Manuscripta Math. \textbf{66} (1990), 309-318.
\bibitem{T}
   K. Tapp, \emph{Conditions of nonnegative curvature on vector bundles and sphere bundles},
   Duke Math. J., to appear.
\bibitem{CodimTwo}
  G. Walschap, \emph{Nonnegatively curved manifolds with souls of codimension 2},
  J. Diff. Geom. \textbf{27} (1988), 525-537.
\end{thebibliography}

\end{document}